\documentclass[12pt]{amsart}

\usepackage{xcolor}
\usepackage[normalem]{ulem}
\usepackage[T1]{fontenc}

\usepackage[margin=2.5cm]{geometry}
\usepackage{amsthm}
\usepackage{amsmath}
\usepackage{hyperref}
\usepackage{dsfont} 

\theoremstyle{definition}
\newtheorem{definition}{Definition}
\newtheorem{theorem}[definition]{Theorem}

\newtheorem{proposition}[definition]{Proposition}

\newcommand{\am}{\textrm{am}}
\newcommand{\sn}{\textrm{sn}}
\newcommand{\cn}{\textrm{cn}}
\newcommand{\dn}{\textrm{dn}}
\newcommand{\strS}{\mathcal{S}}
\newcommand{\A}{\hat{A}}
\newcommand{\B}{\hat{B}}
\newcommand{\D}{\hat{D}}
\newcommand{\E}{\hat{E}}
\newcommand{\Ddd}{\cn^2\big(z^k_q, e \big)+b^2\sn^2\big(z^k_q, e \big)}
\newcommand{\Aaa}{\sn^2\big(z^k_q, e \big)+b^2\cn^2\big(z^k_q, e \big)}
\newcommand{\baba}{\sn^2\big(z^k_q, e \big)\cn^2\big(z^k_q, e \big)}

\newcommand{\Dddx}{\cn^2\big(\xi(t), e \big)+b^2\sn^2\big(\xi(t), e \big)}
\newcommand{\Aaax}{\sn^2\big(\xi(t), e \big)+b^2\cn^2\big(\xi(t), e \big)}
\newcommand{\babax}{\sn^2\big(\xi(t), e \big)\cn^2\big(\xi(t), e \big)}

\title[Lazutkin coordinates of...]{Lazutkin coordinates of the maximal symmetric periodic orbits on the ellipse}
\author{Klaudiusz Czudek}
\address{Klaudiusz Czudek, Institute of Science and Technology Austria (ISTA), Am Campus 1, 3400 Klosterneuburg}
\address{ Institute of Applied Mathematics, Faculty of Physics and Applied Mathematics, Gda{\'n}sk University of Technology, ul. Gabriela Narutowicza 11/12, 80-223 Gda{\'n}sk, Poland}
\email{klaudiusz.czudek@gmail.com}



\begin{document}

\begin{abstract}
In De Simoi J., Kaloshin V., Wei Q. “Dynamical spectral rigidity among $\mathbb{Z}_2$-symmetric strictly convex domains close to a circle (Appendix B coauthored with H. Hezari)”. Ann. of
Math. 186.1 (2017), pp. 277--314  deformational spectral rigidity of $\mathbb{Z}_2$ symmetric domains close to the circle has been shown. One of the steps of the proof was to express the maximal symmetric periodic orbits in Lazutkin parametrization. Here using the action-angle variables we find the second order approximation of Lazutkin coordinates of the maximal symmetric periodic orbits on the ellipses.
\end{abstract}

\maketitle

\section{Introduction}

Let $\Omega$ be a $C^1$ domain in the plane, and let us consider billiard inside $\Omega$, i.e. a particle moving along straight lines that bounces off the boundary with the rule the angle of incidence equals the angle of reflection. Having that the meaning of \textit{periodic orbit} is clear. \textit{The rotation number} of a periodic orbit equals $p/q$, where $q$ is its period and $p$ is the winding number, i.e. the number of times that the associated polygon wraps around the boundary $\partial \Omega$. \textit{The length spectrum} is defined as
$$
\mathcal{L}(\Omega) = \mathbb{N}\{\textrm{length of all periodic orbits in $\Omega$}\}\cup \mathbb{N} \{|\partial \Omega|\}.
$$
A one-parameter family $(\Omega_t)_{t\in[-1,1]}$ of $C^1$ domains is called \textit{	isospectral} if $\mathcal{L}(\Omega_t)=\mathcal{L}(\Omega_0)$ for every $t\in [-1,1]$. A domain $\Omega$ is \textit{dynamically spectrally rigid} if every isospectral one-parameter family of domains with $\Omega_0=\Omega$ is trivial (i.e. $\Omega_t$ is $\Omega_0$ up to isometry).

It is an important open question if every sufficiently smooth domain $\Omega$ is dynamically spectrally rigid. Important progress towards the solution has been made by De Simoi, Kaloshin, Wei \cite{DKW} where the authors show that every sufficiently smooth domain with axial symmetry, sufficiently close to the circle, is dynamically spectrally rigid if we admit only one-parameter families with axial symmetry (see also the introduction of \cite{DKW} for more detailed exposition of known results).

The proof of the main result of \cite{DKW} doesn't incorporate all the periodic orbits in $\Omega$ but just some orbits of a specific type. Namely, for any axially symmetric $\Omega$ and $q\ge 2$ there exists a symmetric periodic orbit of rotation number $1/q$ passing through the marked point of $\partial \Omega$ and having maximal length among other symmetric periodic orbits of rotation number $1/q$ passing through the marked point. We call it the \textit{maximal symmetric $q$-periodic orbit} (Lemma 3.4 \cite{DKW}).

Another important ingredient of \cite{DKW} is the Lazutkin parametrization. Let $\Omega$ be a $C^1$ domain. The Lazutkin parametrization of $\Omega$ is the parametrization of the form
\begin{equation}
\label{E:Lazutkin_general}
x(s) =C_\Pi \int_0^s \rho(s')^{-2/3}ds', \quad C_\Pi =\bigg[ \int_{\partial \Omega} \rho(s')^{-2/3}ds'\bigg]^{-1},
\end{equation}
where $\rho$ is the radius of curvature and $s$ is the arc-length parametrization. Lazutkin weight is a function
$$
\mu(x)=\frac{1}{2 C_\Pi \rho(x)^{1/3}}.
$$

It has been proven (Lemma 5.1 in \cite{DKW}) that for every $\Omega\in \strS^r$, $r\ge 8$, there exist $C^{r-4}$ functions $\alpha$ (odd) and $\beta$ (even) such that for any marked symmetric maximal $q$-periodic orbit $(x_q^0, \cdots, x_q^{q-1})$
\begin{equation}
\label{E:def_alpha}
x^k_q =  \frac{k}{q} + \frac{\alpha(k/q)}{q^2}+ O(q^{-4}).
\end{equation}
Moreover, if $\vartheta_q^k$ represents the angle of reflection of the trajectory at the $k$-th collision then
\begin{equation}
\label{E:def_beta}
\vartheta_q^k = \frac{\mu(x^k_q)}{q} \bigg( 1+ \frac{\beta(k/q)}{q^2}+ O(q^{-4}) \bigg).
\end{equation}

\begin{theorem}
\label{T}
Let us consider the ellipse
\begin{equation}
\label{E:ellipse_def}
\mathcal{E}_e = \bigg\{ (x,y)\in \mathbb{R}^2 : x^2+\frac{y^2}{b^2} =1\bigg\},
\end{equation}
where $e=\sqrt{1-b^2}\in (0,1)$ represents the eccentricity. Then on $\Omega=\mathcal{E}_e$ the functions $\alpha$ and $\beta$ defined by \eqref{E:def_alpha} and \eqref{E:def_beta} take the forms
$$
\alpha(t)=\frac{-eK(e)(1-e^2)}{2\sqrt{1-e^2 \sn^2(K(e)(1-4t), e))}}\bigg[(1-4t)K'(e)\dn(K(e)(1-4t), e)+\partial_2\am(K(e)(1-4t), e) \bigg].
$$
and 
\begin{multline*}
\beta(t)=\alpha'(t) - \frac{8K(e)^2 e^2}{3} \bigg(  \cn^2 \big(\xi(t), e \big)
-\sn^2 \big(\xi(t), e \big)  \bigg) 
+\frac{16K(e)^2e^4 \sn^2 \big(\xi(t), e \big) \cn^2\big(\xi(t), e \big)}{\dn^2 \big(\xi(t), e \big) }\\
+\frac{16K(e)}{1-e^2} \E_1(t)\dn^5\big(\xi(t), e \big),
\end{multline*}
where
\begin{multline*}
\hat{E}_1(t)=
\frac{1}{24} \frac{e^4\sn^2\big(\xi(t), e \big)\cn^2\big(\xi(t), e \big)}
{[\cn^2\big(\xi(t), e \big)+b^2\sn^2\big(\xi(t), e \big)]^2}\\
-\frac{1}{24}\frac{\Aaax}{\Dddx}\\
-\bigg[ 
e^4\frac{\babax}{\big( \Dddx \big)^2} - \frac{1}{3}
\bigg]\\
\cdot
\bigg[ 
\frac{e^4 \babax}{\big( \Dddx \big)^2}
-\frac{\Aaax}{4\big( \Dddx \big)} + \frac{1}{3}
\bigg]\\
-\bigg[
\frac{\Aaax}{4\big( \Dddx \big)} - \frac{1}{3}
\bigg]^2\\
+e^4\frac{\babax}{\Dddx}
\bigg[ \frac{\Aaax}{4 \babax } - \frac{1}{3}
\bigg].
\end{multline*}
Here $K$ is the complete elliptic integral of the first kind, $\dn$ is the Jacobi delta function, $\am$ is the Jacobi amplitude\footnote{The definitions are given in Section 2.}, $\partial_2$ represents the derivative with respect to the second variable and $\xi(t)=K(e)(1-4t)$.
\end{theorem}

\section{Notation and basic facts}

The elliptic parametrization of \eqref{E:ellipse_def} is given by
\begin{equation}
\label{E:def_elliptic}
\gamma(\varphi)=(\sin \varphi, b\cos \varphi), \varphi \in [0,2\pi).
\end{equation}
The family o confocal ellipses $C_\lambda$, $\lambda\in (0, b)$, is given by
\begin{equation}
\label{E:def_confocal}
C_\lambda = \bigg\{ (x,y)\in\mathbb{R}^2 : \frac{x^2}{1-\lambda^2}+\frac{y^2}{b^2-\lambda^2} \bigg\}.
\end{equation}

By $F(\varphi, k)$ and $K(k)$ we denote, respectively, the incomplete and complete Jacobi elliptic integral of the first kind
$$
F(\varphi, k) = \int_0^\varphi \frac{d\tau}{\sqrt{1-k^2\sin^2\tau}}, \quad K(k)=F\bigg(\frac{\pi}{2}, k\bigg).
$$
Jacobi amplitude is defined by the relation $\am(u, k) = \varphi$ if and only if $F(\varphi, k)=u$. Further, the elliptic sine, cosine and delta amplitude are defined as
$$
\sn(u, k) = \sin \am(u,k), \quad \cn(u, k)=\cos\am(u,k), \quad \dn(u,k)= \partial_u \am(u, k).
$$
Let
\begin{equation}
\label{E:def_action_angle}
k_\lambda^2=\frac{e^2}{1-\lambda^2}, \quad \delta_\lambda=  F(\arcsin(\lambda/ b), k_\lambda), \quad \omega_\lambda=\frac{\delta_\lambda}{2K(k_\lambda)}.
\end{equation}

Let $P_1, P_2, \cdots$ be the collision points(with the boundary $\mathcal{E}_e$) of a certain billiard trajectory inside $\mathcal{E}_e$. Let us denote by $P_mP_{m+1}$ the segment joining $P_m$ and $P_{m+1}$. It is a well-known fact (see e.g. \cite{CFS_82}) that if $P_1P_2$ is tangent to one of the confocal ellipses $C_\lambda$  then $P_mP_{m+1}$ is tangent to $C_\lambda$ for every  $m\ge 1$. The next Proposition introduces a parametrization that is useful to study billiard trajectories whose all segments are tangent to $C_\lambda$.

\begin{proposition}[\cite{Chang_Friedberg_88}]
\label{P:action_angle}
Let $\lambda\in (0, b)$ and let us denote
$$S_\lambda (\varphi) = \frac{1}{4}-\frac{F(\frac{\pi}{2}-\varphi, k_\lambda)}{4K(k_\lambda)}.$$
Let $\varphi_0, \varphi_1, \cdots$ be the collision points of a certain billiard trajectory inside $\mathcal{E}_e$ whose all segments are tangent to $C_\lambda$. Then $S_\lambda(\varphi_k)=S_\lambda(\varphi_{0})+k\omega_\lambda$ for every $k\ge 0$.
\end{proposition}

One the ellipse it is convenient to express the Lazutkin parametrization in terms of the elliptic parametrization on $\mathcal{E}_e$. Then it takes the form
\begin{equation}
\label{E:Lazutkin_ellipse}
x(\varphi)=\frac{1}{4} - \frac{F(\frac{\pi}{2}- \varphi, e)}{4K(e)}, \quad x(\varphi)\in [0,1).
\end{equation}
Observe \eqref{E:Lazutkin_ellipse} is the limit of the action angle variables $S_\lambda(\varphi)$ as $\lambda\to 0$. The Lazutkin weight takes the form (by 16.9.1 \cite{Abramowitz_Stegun})
\begin{equation}
\label{E:lazutkin_weight}
\mu(x) = \frac{2K(e)\sqrt{1-e^2}}{\big( \cn^2(K(e)(1-4x), e)+b^2\sn^2(K(e)(1-4x, e)) \big)^{1/2}}= \frac{2K(e)\sqrt{1-e^2}}{|\dn(K(e)(1-4x, e))| }.
\end{equation}

As explained in Introduction we are going to find the Lazutkin coordinates of orbits with the rotation number of the form $1/q$. Thus it is reasonable to define $\lambda_q$ as the unique number with $\omega_{\lambda_q} = \frac{1}{q}$. Then by Proposition \ref{P:action_angle}
\begin{equation}
\label{E:lambdaq}
S_{\lambda_q} (\varphi_q^k) = \frac{k}{q},
\end{equation}
 where $\varphi^k_q$ is the elliptic coordinate of the $k$-th collision point if the maximal symmetric periodic orbit with the rotation number $1/q$.

\section{Calculation of $\alpha$}
\subsection{The Taylor expansion of $\omega_\lambda$ and $\lambda_q$}
Although it is we do not stress that in the notation many of the expressions below depend on the eccentricity $e$ that should be treated as a fixed parameter. Taylor expansion of $F$ at $\varphi=0$, $k=e$ is
\begin{equation}
\label{E:Fat0}
F(\varphi, k) = \varphi + \frac{1}{6}\varphi^3 e^2+\frac{1}{3}e\varphi^3(k-e)+O(\|(\varphi, k-e)\|^5).
\end{equation}
We need also the expansions
\begin{equation}
\label{E:arcsin}
\arcsin\bigg( \frac{\lambda}{\sqrt{1-e^2}} \bigg) = \frac{\lambda}{(1-e^2)^{1/2}}+\frac{\lambda^3}{6(1-e^2)^{3/2}}+O(\lambda^5)) \quad \textrm{at $\lambda=0$}
\end{equation}
and
\begin{equation}
\label{E:modulus}
k_\lambda = \frac{e}{(1-\lambda^2)^{1/2}} = e\bigg( 1+ \frac{\lambda^2}{2}+ O(\lambda^4)\bigg) \quad \textrm{at $\lambda=0$.}
\end{equation}
Plugging \eqref{E:modulus} into the Taylor expansion of $K(k)$ at $k=e$ gives
\begin{equation}
\label{E:modulus2}
K(k_\lambda) = K(e) + \frac{e}{2}K'(e) \lambda^2 + O(\lambda^4) \quad \textrm{at $\lambda=0$.} 
\end{equation}
The Taylor series
$$
\frac{1}{1+x}= 1-x+x^2-x^3+O(x^4)
$$
yields
\begin{equation}
\label{E:modulus3}
\frac{1}{2K(k_\lambda)} = \frac{1}{2K(e)} - \frac{eK'(e)}{4K(e)^2} \lambda^2 + O(\lambda^4) \quad \textrm{at $\lambda=0$.}
\end{equation}
Plugging \eqref{E:arcsin} and \eqref{E:modulus} to \eqref{E:Fat0} we get (see \eqref{E:def_action_angle})
\begin{equation}
\label{E:Fat1}
\delta_\lambda = \frac{\lambda}{(1-e^2)^{1/2}}+\frac{\lambda^3(e^2+1)}{6(1-e^2)^{3/2}} + O(\lambda^5) \quad \textrm{at $\lambda=0$.}
\end{equation}
Multiplying \eqref{E:Fat1} by \eqref{E:modulus3} we arrive at
\begin{equation}
\label{E:omega}
\omega_\lambda = \frac{\lambda}{2K(e)\sqrt{1-e^2}} \bigg[ 1 + \lambda^2 \bigg( \frac{-eK'(e)}{2K(e)}+\frac{e^2+1}{6(1-e^2)} \bigg) + O(\lambda^4) \bigg] \quad \textrm{at $\lambda=0$.}
\end{equation}
The application of the formula 3.6.25 \cite{Abramowitz_Stegun} (p. 16) yields
\begin{equation}
\label{E:omega2}
\lambda_q = \frac{2 K(e) \sqrt{1-e^2}}{q}\bigg( 1+ \frac{1}{q^2} \big( 2K(e)\sqrt{1-e^2} \big)^2 \bigg( \frac{eK'(e)}{2K(e)}-\frac{e^2+1}{6(1-e^2)} \bigg) +O(q^{-4}) \bigg)
\end{equation}

\subsection{Elliptic coordinates of the maximal symmetric periodic orbit}

Inverting the formula $\varphi\to S_\lambda(\varphi)$ from Proposition \ref{P:action_angle} yields
$$\varphi^k_q = \am \bigg( K(k_{\lambda_q}) \bigg( 1- \frac{4k}{q} \bigg), k_{\lambda_q} \bigg).$$
If $q\to \infty$, then $\lambda_q\to 0$ and $k_{\lambda_q}\to e$. Therefore using the Taylor expansion of Jacobi amplitude and \eqref{E:omega2} gives
\begin{multline*}
\label{E:elliptic_periodic}
\varphi^k_q= \am \bigg( K(e) \bigg( 1- \frac{4k}{q} \bigg), e  \bigg)\\
+\frac{1}{q^2} \bigg[ 
2eK'(e)K(e)^2(1-e^2)\dn\bigg( K(e) \bigg( 1- \frac{4k}{q} \bigg), e  \bigg) \bigg(1-\frac{4k}{q} \bigg)\\
+
2eK(e^2)(1-e^2)\partial_2 \am\bigg( K(e) \bigg( 1- \frac{4k}{q} \bigg), e  \bigg)
\bigg]+O(q^{-4}),
\end{multline*}
where $\partial_2$ stands for the partial derivative with respect to the second variable.

\subsection{The Lazutkin coordinates of the maximal symmetric periodic orbit}

By \eqref{E:Lazutkin_ellipse}
$$
x^k_q = \frac{1}{4}- \frac{F(\frac{\pi}{2}-\varphi^k_q, e)}{4K(e)}
$$
Using again the Taylor expansion of $\varphi\longmapsto F(\varphi, e)$ at $\varphi=\am \bigg( K(e) \bigg( 1- \frac{4k}{q} \bigg), e  \bigg)$
$$
\partial_1F|_{(\am(K(e)(1-4t), e)}
=\frac{1}{\sqrt{1-e^2 \sn^2(K(e)(1-4t), e))}}, \quad t\in [0,1],
$$
gives
\begin{equation}
\label{E:Lazutkin_final}
x^k_q = \frac{k}{q}+\frac{\alpha(k/q)}{q^2}+O(q^{-4}),
\end{equation}
where
\begin{equation}
\label{E:alpha_final}
\alpha(t)=\frac{-eK(e)(1-e^2)}{2\sqrt{1-e^2 \sn^2(K(e)(1-4t), e))}}\bigg[(1-4t)K'(e)\dn(K(e)(1-4t), e)+\partial_2\am(K(e)(1-4t), e) \bigg].
\end{equation}

\section{The calculation of $\beta$}

In order to find $\beta$ we the expansion of $\frac{\sin \vartheta^k_q}{\mu(x^k_q)}$ in terms of $1/q$, $q\to \infty$, and then use the relation (B.1) from Appendix B \cite{DKW}. To this end let us denote
$$
L(\varphi, \varphi') = \| \gamma(\varphi) - \gamma (\varphi') \|,
$$
where $\|\cdot\|$ is the Euclidean distance on $\mathbb{R}^2$, $\gamma$ is from \eqref{E:def_elliptic}. The function $L$ is the generating function of the billiard map inside $\mathcal{E}_e$, thus
$$
\sin \vartheta_q^k 
= \sqrt{1-\bigg[\frac{\partial_\varphi L(\varphi, \varphi^{k+1}_q)|_{\varphi=\varphi^k_q}}{|\gamma'(\varphi^k_q)|}\bigg]^2}
$$
\begin{equation}
\label{E:initial_beta}
=\sqrt{1-\frac{\big[ (\cos \varphi^{k+1}_q - \cos \varphi^k_q) \sin \varphi^k_q - b^2  (\sin \varphi^{k+1}_q - \sin \varphi^k_q) \cos \varphi^k_q \big]^2}
{ (\sin^2\varphi^k_q + b^2 \cos^2 \varphi^k_q)\cdot \big( (\cos \varphi^{k+1}_q - \cos\varphi^{k}_q   )^2 + b^2( \sin\varphi^{k+1}_q - \sin\varphi^{k}_q  )^2   \big)  }}
\end{equation}

Let
\begin{equation}
\label{E:def_td}
t:= \varphi_q^k, \quad d := \varphi^{k+1}_q - \varphi^k_q.
\end{equation}
Then
\begin{equation}
\label{E:sin_exp}
\sin \varphi^{k+1}_q- \sin \varphi^{k}_q = d \cos t - \frac{d^2}{2}\sin t - \frac{d^3}{6} \cos t + \frac{d^4}{24} \sin t + \frac{d^5}{120} \cos t - \frac{d^6}{720} \sin t + O(d^7),
\end{equation}
\begin{equation}
\label{E:cos_exp}
\cos \varphi^{k+1}_q- \cos \varphi^{k}_q = - d \sin t - \frac{d^2}{2}\cos t + \frac{d^3}{6} \sin t + \frac{d^4}{24} \cos t - \frac{d^5}{120} \sin t - \frac{d^6}{720} \cos t + O(d^7),
\end{equation}
Let
\begin{equation}
\label{E:notation_beta}
\A = \sin^2 t + b^2 \cos^2 t, \B=\sin t \cos t, \D = \cos^2 t + b^2\sin^2 t.
\end{equation}
Then the denominator in the fraction appearing in \eqref{E:initial_beta} is
$$
d^2 \A^2+ d^3 \A \B e^2 + d^4 \bigg( \frac{\A \D}{4} - \frac{\A^2}{3} \bigg) 
- \frac{d^5}{4}  \A \B e^2
+ d^6 \bigg( \frac{2}{45} \A^2 - \frac{1}{24} \A \D \bigg)
+O(d^7).
$$
In a similar manner the numerator is
$$
d^2 \A^2 + d^3 \A \B e^2 + d^4 \bigg( \frac{1}{4} \B^2 e^4 - \frac{1}{3} \A^2 \bigg)
- \frac{d^5}{4}  \A \B e^2
+ d^6 \bigg( \frac{2}{45} \A^2 - \frac{1}{24} \B^2 e^4 \bigg)
+ O(d^7).
$$
After dividing the numerator and denominator by $d^2\A^2$ and putting
\begin{equation}
\label{E:MN_notation}
M_1 =\frac{\B}{\A} e^2, M_2= \frac{\D}{4\A}-\frac{1}{3}, M_3 = \frac{2}{45} - \frac{1}{24} \frac{\D}{\A}, N_2= \frac{\B^2}{4\A^2}e^4- \frac{1}{3}, N_3 = \frac{2}{45} - \frac{1}{24} \frac{\B^2}{\A^2} e^4
\end{equation}
the denominator can be written down as
$$
1+dM_1+d^2M_2- \frac{d^3}{4} M_1 + d^4 M_3 + O(d^5),
$$
while the numerator
$$
1+ dM_1 + d^2N_2 - \frac{d^3}{4}M_1+d^4 N_3 + O(d^5).
$$
Using the Taylor expansion of $1/(1+x)$ the inverse of the denominator can be written down as
\begin{multline*}
1- dM_1 + d^2 (M_1^2-M_2) + d^3\bigg( \frac{M_1}{4}+2M_1M_2 - M_1^3 \bigg)
+ d^4\bigg(-M_3+M_2^2 - \frac{M_1^2}{2} - 3M_1^2 M_2 + M_1^4 \bigg)
+ O(d^5).
\end{multline*}
The product of that with the numerator gives the expansion of the fraction appearing under the square root in \eqref{E:initial_beta} in terms of $d$ and therefore of the entire expression under the square root. In order to obtain the expansion of $\sin\vartheta^k_q$ we factor out $d^2(M_2-N_2)$ and use the Taylor expansion of $\sqrt{1+x}$. After that we arrive at 
\begin{equation}
\label{E:sin1}
\sin \vartheta^k_q = \sqrt{M_2-N_2} \bigg( d  + d^2 M_1 + d^3 \frac{\E_1}{M_2-N_2} \bigg) +O(d^4),
\end{equation}
where 
\begin{equation}
\label{E:Edef}
\E_1 =  M_3-N_3 - N_2M_1^2+N_2M_2 - M_2^2+ M_1^2M_2.
\end{equation}
Now, let us express $d=\varphi^{k+1}_q - \varphi^k_q$ in terms of $1/q$. We 
have
$$
x_q^{k+1} = \frac{k}{q}+\frac{1}{q}+\frac{\alpha(k/q)}{q^2}+ \frac{\alpha'(k/q)}{q^3} + O(q^{-4}),
$$
and therefore
$$
x^{k+1}_q - x^k_q = \frac{1}{q}+\frac{\alpha'(k/q)}{q^3} + O(q^{-4}).
$$
Let $y_q^k=K(e)(1-4x^k_q)$ and $z_q^k=K(e)(1-4k/q)$. By inverting \eqref{E:Lazutkin_ellipse} we arrive at
\begin{multline*}
d
=\varphi_q^{k+1} - \varphi_q^k
= 
4K(e)\dn\big( y^k_q, e \big) \bigg[ x^{k+1}_q - x^k_q
+
2 K(e) e^2 (x^{k+1}_q - x^k_q)^2  \frac{ \sn\big(y^k_q, e \big) \cn\big(y^k_q, e \big)} {\dn\big(y^k_q, e \big)} \\
-
\frac{8}{3} K(e)^2 e^2 (x^{k+1}_q - x^k_q)^3 
 \bigg(  \cn^2 \big(y^k_q, e \big)
-\sn^2 \big(y^k_q, e \big)  \bigg)
\bigg]
+
O\big((x^{k+1}_q - x^k_q)^4\big)\\
=
4K(e)\dn\big(y^k_q, e \big) 
\bigg[ \frac{1}{q}
+\frac{2K(e) e^2}{q^2}  \frac{ \sn\big(z^k_q, e \big) \cn\big(z^k_q, e \big)} {\dn\big(z^k_q, e \big)} \\
+\frac{1}{q^3}\bigg(\alpha'(k/q) - \frac{8K(e)^2 e^2}{3} \bigg(  \cn^2 \big(z^k_q, e \big)
-\sn^2 \big(z^k_q, e \big)  \bigg)  \bigg) \bigg] + O(q^{-4})
\end{multline*}
Since
$$
4K(e)\dn\big(y^k_q, e \big)
= 4K(e)\dn\big(z^k_q, e \big) \\
+ \frac{16K(e)^2 e^2 \alpha(k/q)}{q^2} \sn\big(z^k_q, e\big)\cn\big(z^k_q, e\big) + O(q^{-4})
$$
we have
$$
d^2=4K(e)\dn\big(y^k_q, e \big)
\bigg( 
\frac{4K(e)\dn\big(z^k_q, e \big) }{q^2}
+ \frac{16K(e)^2e^2}{q^3}
\sn\big(z^k_q, e \big) \cn\big(z^k_q, e \big) + O(q^{-4})
\bigg)
$$
and
$$
d^3=4K(e)\dn\big(y^k_q, e \big)\bigg( \frac{4K(e)\dn\big(z^k_q, e \big) }{q^3} +O(q^{-4}) \bigg).
$$
Moreover,
$$
M_1=
\frac{e^2 \sn \big(K(e)\big(1-4k/q\big), e \big) \cn\big(K(e)\big(1-4k/q\big), e \big)}{\dn^2 \big(K(e)\big(1-4k/q\big), e \big) }+O(q^{-2}).
$$
By the definition of $M_2$, $N_2$, \eqref{E:Lazutkin_ellipse} and 16.9.1 in \cite{Abramowitz_Stegun} we have
\begin{multline*}
\sqrt{M_2-N_2} = \frac{\sqrt{1-e^2}}{2\big(\sin^2 \varphi^k_q+b^2\cos^2 \varphi^k_q\big)}
=
\frac{\sqrt{1-e^2}}{2\big(\cn^2(K(e)(1-4x^k_q), e)+b^2\sn^2(K(e)(1-4x^k_q), e)\big)} \\
=\frac{\sqrt{1-e^2}}{2\dn^2(K(e)(1-4x^k_q), e) }.
\end{multline*}
Going back to \eqref{E:sin1} we factor out 
$$\frac{4K(e)\dn\big(K(e)\big(1-4x^k_q\big), e \big)}{q}\sqrt{M_2-N_2}=\frac{\mu(x^k_q)}{q}$$
By (B.1) in \cite{DKW} the term that we are looking for equals
\begin{multline*}
\beta(k/q)=\alpha'(k/q) - \frac{8K(e)^2 e^2}{3} \bigg(  \cn^2 \big(K(e)\big(1-4k/q\big), e \big)
-\sn^2 \big(K(e)\big(1-4k/q\big), e \big)  \bigg)  \\
+\frac{16K(e)^2e^4 \sn^2 \big(K(e)\big(1-4k/q\big), e \big) \cn^2\big(K(e)\big(1-4k/q\big), e \big)}{\dn^2 \big(K(e)\big(1-4k/q\big), e \big) }\\
+\frac{16K(e)}{1-e^2} \E_1(k/q)\dn^5\big(K(e)\big(1-4k/q\big), e \big),
\end{multline*}
where
\begin{multline*}
\hat{E}_1=
\frac{1}{24} \frac{e^4\sn^2\big(z^k_q, e \big)\cn^2\big(z^k_q, e \big)}
{[\cn^2\big(z^k_q, e \big)+b^2\sn^2\big(z^k_q, e \big)]^2}\\
-\frac{1}{24}\frac{\Aaa}{\Ddd}\\
-\bigg[ 
e^4\frac{\baba}{\big( \Ddd \big)^2} - \frac{1}{3}
\bigg]\\
\cdot
\bigg[ 
\frac{e^4 \baba}{\big( \Ddd \big)^2}
-\frac{\Aaa}{4\big( \Ddd \big)} + \frac{1}{3}
\bigg]\\
-\bigg[
\frac{\Aaa}{4\big( \Ddd \big)} - \frac{1}{3}
\bigg]^2\\
+e^4\frac{\baba}{\Ddd}
\bigg[ \frac{\Aaa}{4\big( \baba \big)} - \frac{1}{3}
\bigg].
\end{multline*}

\section{Acknowledgments}
The author thanks Vadim Kaloshin, Corentin Fierobe and Jacopo De Simoi for useful discussions.

\bibliographystyle{plain}
\bibliography{Bibliography}
\end{document}